\documentclass{amsart}

\usepackage{amsmath}
\usepackage{graphicx}
\usepackage{url}
\usepackage{mathrsfs}

\newtheorem {theo} {Theorem}
\newtheorem {prop} {Proposition}
\newtheorem {defi} {Definition}
\newtheorem {conj} {Conjecture}

\DeclareMathAlphabet{\nim}{U}{bbold}{m}{n}

\title{Analysis of mis\`ere Sprouts game with reduced canonical trees}
\author{Julien Lemoine - Simon Viennot}

\begin{document}

\begin{abstract}
Sprouts is a two-player topological game, invented in 1967 by Michael Paterson and John Conway. The game starts with $p$ spots, lasts at most $3p-1$ moves, and the player who makes the last move wins. In the  \textit{mis\`ere} version of Sprouts, on the contrary, the player who makes the last move \textit{loses}.\par
Sprouts is a very intricate game, and the first computer analysis in 1991 reached only $p=11$. New results were made possible in 2007 up to $p=32$ by using combinatorial game theory: when a position is a sum of independant games, it is possible to replace some of these games by a natural number, called the \textit{nimber}, without changing the winning or losing outcome of the complete position.\par
However, this reduction does not apply to the mis\`ere version, making the analysis of Sprouts (and more generally of any game) more difficult in the mis\`ere version. In 1991, only $p=9$ was reached in mis\`ere Sprouts, and we describe in this paper how we obtained up to $p=17$.\par
First, we describe a theoretical tool, \textit{the reduced canonical tree}, which plays a role similar to the nimber in the normal version. Then, we describe the way we have implemented it in our program, and detail the results it allowed us to obtain on mis\`ere Sprouts.
\end{abstract}

\maketitle

\section{Introduction}

\textit{Sprouts} is a two-player game, which needs only a sheet of paper and a pen to play, with extremely simple rules: for example, Martin Gardner's article of 1967 \cite{mg67} (when the game was invented) is a good introduction to start the study of the game.\par

Sprouts is a \textit{combinatorial game}: two players play alternately, knowing all the possible information to choose their next move. There is no room for chance in the game. Moreover, Sprouts is an \textit{impartial} combinatorial game: from any position, the same moves are available to either player.\par

In the \textit{normal} version of the game, the winner is determined by the following rule: a player who cannot make a move loses. Draws are not possible, and since a game beginning with $p$ spots is finite (with at most $3p-1$ moves), there must be a winning strategy for one of the player (but for that, of course, he needs to play perfectly).\par

\begin{defi}
$S_p^+$ denotes the normal version of the Sprouts game starting with $p$ spots (and $S_p^-$ the mis\`ere version).
\end{defi}

Finding which player has a winning strategy is difficult because of the game complexity. The first manual analysis only achieved to solve $S_6^+$, and it required to consider a lot of cases through many pages of reasoning. In 1991, the first program of Sprouts enabled Applegate, Jacobson and Sleator \cite{ajs91} to extend this analysis up to $S_{11}^+$, and to formulate the following conjecture:

\begin{conj}
The first player has a winning strategy in $S_p^+$ if and only if $p$ is $3$, $4$ or $5$ modulo $6$.
\end{conj}

Later computation \cite{lv07} in 2007 proved this conjecture to be true up to $S_{32}^+$.\par

In the \textit{mis\`ere} version of the game, the last player able to move is this time the loser. Games in mis\`ere version are almost always more difficult to analyze. For example, in the normal version, a disjunctive sum of two losing positions is losing. But in the mis\`ere version, the sum can be losing or winning, depending on the case. Detailed explanations about the difficulty of mis\`ere games are given in the \textit{Theory} section.

Because of this difficulty, Applegate et al. in 1991 only reached $S_9^-$ in the mis\`ere version, upon which they formulated the following conjecture:

\begin{conj}
(false) The first player has a winning strategy in $S_p^-$ if and only if $p$ is equal to $0$ or $1$ modulo $5$.
\end{conj}

This conjecture was invalidated by Josh Purinton (unpublished work), who computed first that the winning strategy for $S_{12}^-$ is actually for the first player\footnote{see http://www.wgosa.org/confirmation12.htm -- Josh Purinton has also solved the mis\`ere version up to $S_{16}^-$.}. In fact, we observe that above $S_{5}^-$, we get back to a pattern of period 6, just as in the normal version, but with a shift: known values up to now are shown in the table of figure \ref{WL17}.

\begin{figure}[h]
\centering
\begin{tabular}{ |*{18}{c|} }
\hline
$p$ & 1 & 2 & 3 & 4 & 5 & 6 & 7 & 8 & 9 & 10 & 11 & 12 & 13 & 14 & 15 & 16 & 17\tabularnewline
\hline
$S_p^-$ & W & L & L & L & W & W & L & L & L & W & W & W & L & L & L & W & W\tabularnewline
\hline
\end{tabular}
\caption{Computed Win/Loss state of $S_p^-$.}
\label{WL17}
\end{figure}

We can formulate the following conjecture:

\begin{conj}
The first player has a winning strategy in $S_p^-$ if and only if $p$ is equal to $0$, $4$ or $5$ modulo $6$ -- except if $p=1$ or $4$.
\end{conj}

We describe in this article how we achieved to analyze the game up $S_{17}^-$, beginning first with some well-known notions of combinatorial game theory. Then, we explain how we implemented them in a program, and finally, we will give details on the computed results.\par

In this article, we use the notation of Sprouts positions described in \cite{lv07}, which is derived from the notation of \cite{ajs91}.

\section{Theory}

In this section, we give an overview of some well known results of the theory of mis\`ere games, without proving again all of them. A good entry point for this theory is \textit{On Numbers And Games} \cite{jc01}, a book from Conway (chapter 12, pp. 136--152), or the famous \textit{Winning Ways} \cite{ww01} from Berkelamp, Conway and Guy (chapter 13, pp. 413--452).

\subsection{Indistinguishability}

We give first some standard definitions from the theory of impartial combinatorial games. A game in a given state will be called a \textit{position}. The \textit{outcome} of a position is W (win) or L (loss), depending on the existence of a winning strategy from this position. The \textit{sum} of two positions $\mathscr{P}_1$ and $\mathscr{P}_2$ is the position obtained by grouping them, so that at each move, the player can choose to move in the first or the second one. This sum will be denoted $\mathscr{P}_1+\mathscr{P}_2$.\par

Following \cite{tp05}, we say that two positions $\mathscr{P}_1$ and $\mathscr{P}_2$ are \textit{indistinguishable}, and we note $\mathscr{P}_1\sim \mathscr{P}_2$, if for any other position $\mathscr{T}$, the positions $\mathscr{P}_1+\mathscr{T}$ and $\mathscr{P}_2+\mathscr{T}$ have the same outcome. This notion is of practical interest because if two Sprouts positions are known to be indistinguishable, it is possible to replace the biggest one by the smallest one in order to accelerate computation.

Indistinguishability depends on the considered game, and also on its version: two positions can be indistinguishable in the normal version, but not in the mis\`ere one. For example, in the case of Sprouts, the start positions with 1 and 2 spots $S_1$ and $S_2$ are indistinguishable in the normal version (let us denote $S_1\sim_+S_2$). But in the mis\`ere version the empty position $\mathscr{T}$ distinguish them, since $S_1$ is a Win and $S_2$ a Loss ($S_1\not\sim_-S_2$).

\subsection{Indistinguishability in the normal version}

Indistinguishability allows us to simplify considerably any impartial game in the normal version (when the player who cannot move is the loser). But, first, we need to go back to a very simple and essential game: the game of Nim.

\begin{defi}
A position $\mathscr{P}$ is entirely defined by the set of possible moves from it. We note $\mathscr{P}=\{\mathscr{P}_1,\mathscr{P}_2,\mathscr{P}_3,...\}$, where $\mathscr{P}_i$ are the children of $\mathscr{P}$.
\end{defi}

\begin{defi}
 Let $\nim{n}$ denote the Nim-heap with $n$ matchsticks\footnote{Children should never play with matchsticks or any other source of fire.}.
\end{defi}

\begin{itemize}
 \item $\nim{0}$ is a void heap. This is a terminal position, where no move is possible: $\nim{0}=\{\}$.
 \item $\nim{1}$ is a heap of only one matchstick. The only possible move is to remove this matchstick, so $\nim{0}$ is the only option: $\nim{1}=\{\nim{0}\}$.
 \item $\nim{2}$ is a heap of two matchsticks: it is possible to remove one or the two matchsticks, and $\nim{2}=\{\nim{0};\nim{1}\}$.
 \item The general rule, for any positive number $n$ is $\nim{n}=\{\nim{0};\nim{1};...;\nim{n-1}\}$.
\end{itemize}

The importance of Nim comes from the following result:

\begin{theo}
(of Sprague-Grundy) For any impartial combinatorial game, all the positions are indistinguishable from some Nim-heap, called the \emph{nimber} of the position.
\end{theo}

In the case of Sprouts, $S_2\sim_+ \nim{0}$ and \verb/ABCD.}AB.}CD.}]!/ $\sim_+ \nim{3}$, which means that the nimber of $S_2$ is $\nim{0}$ and the one of \verb/ABCD.}AB.}CD.}]!/ is $\nim{3}$.

Indistinguishability is an equivalence relation, and the corresponding equivalence classes are called \textit{indistinguishability classes}. The Sprague-Grundy theorem then states that, in the normal version of impartial games, there are very few and particularly simple indistinguishability classes. This greatly simplifies the analysis of games where positions appear as the sum of smaller ones.

Unfortunately, this theorem does not apply to the mis\`ere version. There are many more indistinguishability classes, and John Conway proves in \cite{jc01} that instead of the nimber, we need the concept of \textit{reduced canonical tree} to analyze the classes of the mis\`ere version. We describe it in the following sections.

\subsection{Game tree}

We call \textit{game tree obtained from a position $\mathscr{P}$} the tree where the vertices are all the positions that can be reached by playing moves from $\mathscr{P}$, and where two positions $\mathscr{P}_1$ and $\mathscr{P}_2$ are linked by an edge if $\mathscr{P}_2$ is obtained from $\mathscr{P}_1$ in only one move.

In order to construct the game tree $\mathscr{A}$ obtained from a position $\mathscr{P}$, we need of course to know the rules of the game, so that we can compute the children of $\mathscr{P}$, but it is important to note that it is then possible to study directly the game tree $\mathscr{A}$, without refering anymore to the underlying position, or even the underlying game.

\subsection{Canonical trees}

When two branches of a game tree are perfectly identical, it means that the player can choose between two moves leading exactly to the same situation. Choosing one move or the other will not make any difference in the game, so redundant branches in a game tree are useless. We call \textit{canonical tree} the game tree where all redundant branches have been deleted.

\begin{figure}[ht]
\centering
\includegraphics[scale=0.4,angle=-90]{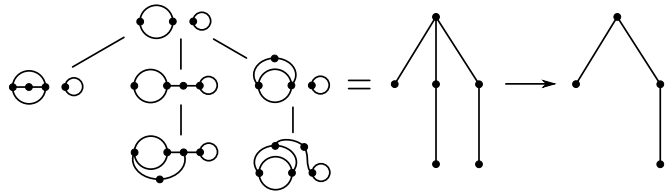}
\caption{Game tree obtained from a given position and corresponding canonical tree}
\label{trees}
\end{figure}

Figure \ref{trees} shows on the left a game tree obtained from a given Sprouts position. The two branches on the right lead to similar games, so we can merge them into a single branch to obtain the canonical tree (on the right). The canonical tree, as well as the game tree, is of height 2, because the longest possible game ends in 2 moves.

The canonical tree corresponding to a given game tree $\mathscr{A}$ can be defined recursively:
\begin{itemize}
 \item compute the canonical tree of each child of $\mathscr{A}$.
 \item in these canonical children, delete all the redundant ones.
\end{itemize}

\begin{figure}[ht]
\centering
\includegraphics[scale=0.4,angle=-90]{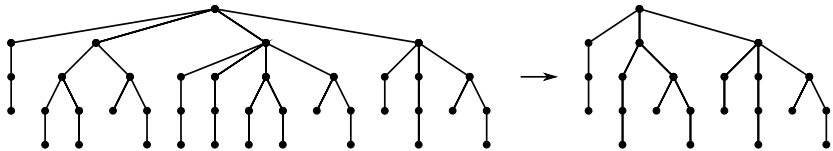}
\caption{Canonization of a game tree}
\label{CTp2}
\end{figure}

This notion of canonical tree allows us to keep only the necessary and sufficient information of a game tree needed to describe a position: if two positions have the same canonical tree, then they are indistinguishable whatever the version of the game is (normal, mis\`ere, or any other rule where the result depends only on the number of played moves).

Figure \ref{CTp2} shows the canonization of the game tree resulting from the Sprouts position \verb/0.AB.}AB.}]!/ (this position comes from $S_2$, after linking a spot to itself), obtained with our program.

\subsection{Reversible moves}

\subsubsection{Reversible moves}

\begin{defi}
Let $\mathscr{G}=\{\mathscr{G}_1,\mathscr{G}_2,\mathscr{G}_3,...\}$ be a non-empty canonical tree. Then, we say that a canonical tree $\mathscr{H}$ is obtained from $\mathscr{G}$ by adding reversible moves if $\mathscr{H}=\{\mathscr{G}_1,\mathscr{G}_2,\mathscr{G}_3,...,\mathscr{R}_1,\mathscr{R}_2,\mathscr{R}_3,...\}$ and each tree $\mathscr{R}_j$ has $\mathscr{G}$ in its children (the moves $\mathscr{R}_j$ are said \emph{reversible}).
\end{defi}

$\mathscr{G}$ and $\mathscr{H}$ have the same outcome: if $\mathscr{G}$ is a Win, it means one of its children $\mathscr{G}_i$ is a Loss. Since $\mathscr{G}_i$ is also a child of $\mathscr{H}$, it implies that $\mathscr{H}$ is a Win. Conversely, if $\mathscr{G}$ is a Loss, it means all children $\mathscr{G}_i$ are wins. Moreover, each $\mathscr{R}_i$ is a Win, since it has $\mathscr{G}$ as a child, so all children of $\mathscr{H}$ are wins, and $\mathscr{H}$ is a Loss.

But in fact, there is a stronger result, which states that $\mathscr{G}$ and $\mathscr{H}$ are indistinguishable in the mis\`ere version. Indeed, if a player has a winning strategy for the sum $\mathscr{G}+\mathscr{T}$, then he can win $\mathscr{H}+\mathscr{T}$, by playing the same moves, and extending the strategy to the following case:
\begin{itemize}
\item if the opponent plays one of the move $\mathscr{R}_j$, he should answer by playing the move $\mathscr{G}$ for this component (he ``reverses'' the move $\mathscr{R}_j$, hence the name \textit{reversible move}).
\end{itemize}

\begin{figure}[ht]
\centering
\includegraphics[scale=0.4]{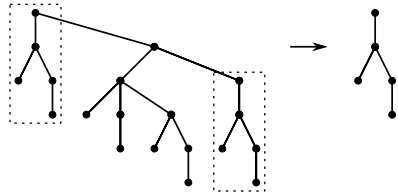}
\caption{Example of reversible move}
\label{reversible}
\end{figure}

Figure \ref{reversible} shows an example of reversible move: $\mathscr{H}$ is the tree on the left, and $\mathscr{G}$ the one on the right. There is a reversible move $\mathscr{R}_1$, from which two moves are possible, and one of them is $\mathscr{G}$.

\subsubsection{Particular case of the empty tree}
\label{emptytree}

If $\mathscr{G}$ is the empty tree, the above definition holds, but we need an additional clause to ensure that  $\mathscr{G}$ and $\mathscr{H}$ are indistinguishable in the mis\`ere version: $\mathscr{H}$ must be a Win in the mis\`ere version.

The reason is that when the two players play the game $\mathscr{H}+\mathscr{T}$ as described above, a particular case arises when $\mathscr{T}$ is finished before playing any move in $\mathscr{H}$. If $\mathscr{G}$ is not empty, the strategy for playing $\mathscr{H}$ is described above.

But if $\mathscr{G}$ is empty, the game on $\mathscr{G}+\mathscr{T}$ is supposed to be finished, and the player to play should be the winner. However, he finds himself forced to play one of the reversible moves $\mathscr{R}_j$, and of course, he can win in this case only if $\mathscr{H}$ is a Win in the mis\`ere version. The clause (also known as the \textit{proviso}) ensures that $\mathscr{G}$ and $\mathscr{H}$ have the same outcome even in this particular case.

\subsection{Reduced canonical trees}
\label{rct}

\subsubsection{Reducers}
\label{reducers}

The previous section shows that if $\mathscr{H}$ is obtained from $\mathscr{G}$ by adding reversible moves, then $\mathscr{G}$ and $\mathscr{H}$ are indistinguishable in the mis\`ere version, but what interests us is the reverse way: given a tree $\mathscr{H}$, we try to simplify it by reducing it to a given tree $\mathscr{G}\subset{\mathscr{H}}$, obtained by pruning reversible moves. This kind of tree $\mathscr{G}$ will be called a \textit{reducer} of $\mathscr{H}$.

Let us remark that $\mathscr{G}$ is not necessarily unique, which seems at first to be a problem when there are several possible reducers. However, the reversible moves $\mathscr{R}_1,\mathscr{R}_2,\mathscr{R}_3,...$ must all include $\mathscr{G}$ strictly, so their height is at least the height of $\mathscr{G}$ plus 1. It follows that:

\begin{itemize}
\item A reducer is always of the form: \{children of $\mathscr{H}$ with a height less than a given value\}.
\item Two reducers are always comparable for the inclusion. If there are two different reducers, the smallest one is included in the biggest one, and moreover it is a reducer of the biggest one.
\item There exists a (unique) smallest reducer.
\end{itemize}

\subsubsection{Reduced canonical trees}

We can then define the \textit{reduced canonical tree} of a canonical tree, by pruning all the possible reversible moves. This reduced canonical tree is constructed recursively:
\begin{itemize}
\item compute the reduced canonical tree of each child.
\item in this reduced canonical children, delete the redundant ones.
\item reduce the resulting tree with the smallest possible reducer.
\end{itemize}

\subsection{Indistinguishability in mis\`ere version}

The main interest of these reduced canonical trees is that if two positions have the same reduced canonical tree, then they are indistinguishable in mis\`ere version. It is natural to ask for the converse.

The converse is true for combinatorial games in normal version: we know that if two positions have a different nimber, then they are distinguishable. Indeed, if $\mathscr{P}_1$ and $\mathscr{P}_2$ have different nimbers, then $\mathscr{P}_1$ distinguish them, because $\mathscr{P}_1+\mathscr{P}_1$ is a Loss, while $\mathscr{P}_2+\mathscr{P}_1$ is a Win.

In the case of a mis\`ere game, a result proved in \cite{jc01} (p. 149) seems at first to answer the question: given two different reduced canonical trees $\mathscr{G}$ and $\mathscr{H}$, there exists a reduced canonical tree $\mathscr{T}$ such that $\mathscr{G}+\mathscr{T}$ and $\mathscr{H}+\mathscr{T}$ have different mis\`ere outcomes, which means that $\mathscr{T}$ distinguish $\mathscr{G}$ and $\mathscr{H}$.

 From the above result, we could conclude that in the mis\`ere version, indistinguishability classes are exactly the reduced canonical trees. However, when studying a game in particular, it is possible that some positions with different reduced canonical trees are in fact indistinguishable. Let us give an example.

We consider the mis\`ere game of Nim, restricted to heaps of size $\leq 2$. Since $\nim{1}+\nim{1}\sim_-\nim{0}$, indistinguishability classes are of the form $n\times\nim{2}$ or $n\times\nim{2}+\nim{1}$ ($n\geq 0$). The possible moves from this kind of positions are as follows:

\begin{itemize}
 \item from $n\times\nim{2}$ ($n\geq 1$), we can move to $(n-1)\times\nim{2}+\nim{1}$ or $(n-1)\times\nim{2}$.
 \item from $n\times\nim{2}+\nim{1}$ ($n\geq 1$), we can move to $n\times\nim{2}$, $(n-1)\times\nim{2}+\nim{1}$ or $(n-1)\times\nim{2}$.
\end{itemize}

It enables us to determine recursively that the only losing positions are $\nim{1}$, and $2n\times\nim{2}$ ($n\geq 1$), and then that the only indistinguishability classes are $\nim{0}$ ; $\nim{1}$ ; $\nim{2}$ ; $\nim{2}+\nim{1}$ ; $\nim{2}+\nim{2}$ ; $\nim{2}+\nim{2}+\nim{1}$. Indeed, when a position includes at least 3 times $\nim{2}$, deleting a pair of $\nim{2}$ does not change the outcome.

It implies that even if the reduced canonical trees of $\nim{2}$ and $\nim{2}+\nim{2}+\nim{2}$ are different, since no reduced canonical tree appearing in the game distinguish them, then they are indistinguishable\footnote{This example is detailed in \cite{si08}.}. Actually, this simplification can happen as soon as all the possible reduced canonical trees do not appear in a given game. The new indistinguishability classes, bigger and less numerous, are the root of Thane Plambeck's work (see for example \cite{tp05}) on \textit{mis\`ere quotients}.

Unfortunately, this theory seems difficult to apply to Sprouts, where a lot of different reduced canonical trees appear. Trying to find less numerous indistinguishability classes may still be an option for Sprouts, but we restrained our work to the analysis of reduced canonical trees.

Going back to our example, $\nim{2}$ and $\nim{2}+\nim{2}+\nim{2}$ are on the contrary distinguishable within Sprouts: the position \verb/ABCD.}ABEF.}CDFE.}]!/, whose reduced canonical tree is $\{\nim{1};\{\nim{2}\}\}$, distinguish them, since $\nim{2}+\{\nim{1};\{\nim{2}\}\}$ is a Win, while $\nim{2}+\nim{2}+\nim{2}+\{\nim{1};\{\nim{2}\}\}$ is a Loss.

\subsection{Count}

It is interesting to count the canonical trees and the reduced canonical trees, in order to evaluate their practical interest.

\begin{figure}[ht]
\centering
\includegraphics[scale=0.4,angle=-90]{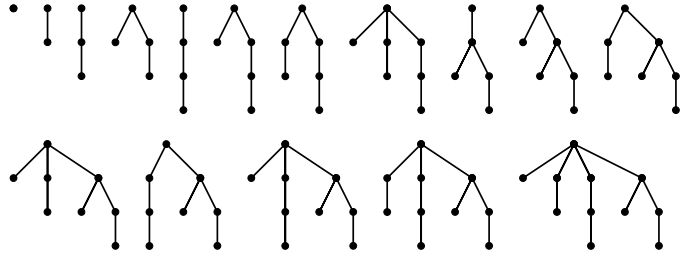}
\caption{Canonical trees of height $\leq 3$}
\label{trees height 3}
\end{figure}

First, we can count the exact number of canonical trees of a given height, as on figure \ref{trees height 3}, which shows the 16 canonical trees of height $\leq 3$.

\begin{prop}
There is $2^{(2^{(...^{(2^0)})})}$ (with $h+1$ times the number 2) canonical trees of height $\leq h$.
\end{prop}

This can be proved by recursion. Let $\mathscr{C}_h$ denote the set of all canonical trees of height $\leq h$ and $c_h$ its cardinal number. The above formula can then be written $c_{h+1}=2^{c_h}$. This comes from the fact that a canonical tree of height $\leq {h+1}$ can be defined as the set of canonical trees of its children, which are of height $\leq h$. It means that there is bijection between the set $\mathscr{C}_{h+1}$ of all canonical trees of height $\leq {h+1}$ and $\mathscr{P}(\mathscr{C}_h)$, the power set of $\mathscr{C}_h$.

It implies that the number of canonical trees of height $\leq n$ is, for increasing $n$: 1;2;4;16;65536;$2^{65536}$... which should be compared with the number of nimbers corresponding to trees of height $\leq n$: 1;2;3;4;5;6...

\begin{figure}[ht]
\centering
\includegraphics[scale=0.4]{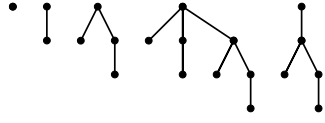}
\caption{Reduced canonical trees of height $\leq 3$: $\nim{0}$ ; $\nim{1}$ ; $\nim{2}$ ; $\nim{3}$ and  $\{\nim{2}\}$}
\end{figure}

Unfortunately, the number of reduced canonical trees is more similar to the first case: 1;2;3;5;22;4171780... (see \cite{gs56}). In fact, there is a lot of reductions for small values, but they rarefy quickly and we come back to a growth of the form $c_{n+1}=2^{c_n}$. For example, $2^{22}=4194304$ is very close to 4171780, and the next term is more than 99,99\% of $2^{4171780}$ (the exact value is given in \cite{jc01} p. 140).

\subsection{Nim-heaps}

The previous paragraph shows that reduced canonical trees appearing in the mis\`ere version of Sprouts are much more numerous and complex than the nimbers of the normal version. However, almost-terminal positions frequently have a very simple reduced canonical tree, namely the tree of a Nim-heap. We detail in the following the properties explaining this fact.

\begin{defi}
The Mex (minimum excluded) of a set of natural numbers is defined as the least natural number which is not included in the set.
\end{defi}

Example: $Mex(0;1;4;3;1;7)=2$.

The canonical trees corresponding to Nim-heaps cannot be simplified with reversible moves, but on the contrary, it is possible to simplify a position whose children are Nim-heaps (this result, as well as the following one, are proved in \cite{jc01} p.139):

\begin{theo}
\label{theomex}
 A position whose children are all Nim-heaps is itself indistinguishable from a Nim-heap, except if all the children have size at least 2. The position is then indistinguishable from the Mex of the children.
\end{theo}

For example, $\mathscr{P}_1=\{\nim{0};\nim{1};\nim{3};\nim{5}\}$ (a position whose children are Nim-heaps of size 0;1;3;5) is indistinguishable from a Nim-heap of size 2: $\mathscr{P}_1=\{\nim{0};\nim{1};\nim{3};\nim{5}\}\sim\nim{2}$.

On the contrary, $\mathscr{P}_2=\{\nim{2};\nim{3}\}$ (a position whose children are Nim-heaps of size 2 and 3) cannot be reduced.

Moreover, Nim-heaps have another interesting property, which means that it is sometimes possible to reduce a position constituted of a sum of Nim-heaps:

\begin{theo}
If at least one of the number $m$ or $n$ is equal to 0 or 1, then $\nim{m}+\nim{n} \sim \nim{q}$, where $q=m\oplus n$.
\end{theo}

``$\oplus$'' is the \textit{Nim-sum}, which is done by a bitwise ``exclusive or'' on the two numbers. For example, $\nim{3}+\nim{1}\sim\nim{2}$, or $\nim{4}+\nim{1}\sim\nim{5}$.

Reduction is not possible if $m$ and $n$ are $\geq 2$. In that case, $\nim{m}+\nim{n}$ is not indistinguishable from a given Nim-heap. For example, $\nim{2}+\nim{2}=\{\nim{2}+\nim{0};\nim{2}+\nim{1}\}\sim\{\nim{2};\nim{3}\}$, and it has already been described above as a position impossible to reduce.

Let us give an example for Sprouts. We consider the game tree obtained from $S_3$. The positions contained in this game tree correspond to 55 different canonical trees, and only 2 of them are not reducible to a Nim-heap:

\begin{itemize}
 \item $S_2=$\verb/0.0.}]!/, which has two children, both indistinguishable from $\nim{2}$. Its reduced canonical tree is then $\{\nim{2}\}$.
 \item \verb/0.0.AB.}AB.}]!/, which is ``contaminated'' by its child $S_2$. Its reduced canonical tree is $\{\nim{1};\{\nim{2}\}\}$.
\end{itemize}

The other 53 positions are all reducible to Nim-heaps, which shows the importance of this concept when analyzing small positions of mis\`ere Sprouts.

\subsection{Recovery with reversible moves}

We have seen in the previous paragraph that the position \verb/0.0.AB.}AB.}]!/ is ``contaminated'' by one of its children. But when going back up the tree, it can happen that some positions are not, i.e. they are indistinguishable from a given Nim-heap, even though some positions of their subtree are not.

For example, this is the case of $S_3$. It has three children, two of them reducible to $\nim{0}$, and the other being \verb/0.0.AB.}AB.}]!/. $S_3$ is then indistinguishable from $\{\nim{0};\{\nim{1};\{\nim{2}\}\}\}$, but this is reducible to $\nim{1}$, because it can be obtained from $\nim{1}$ by adding the reversible move $\{\nim{1};\{\nim{2}\}\}$. It means that $S_3$ is indistinguishable from a Nim-heap, even though it was not the case of one of its children.

\begin{figure}[ht]
\centering
\includegraphics[scale=0.4]{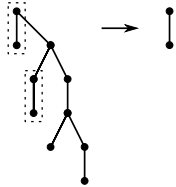}
\caption{Reduced canonical tree obtained from $S_3$}
\end{figure}

This property of ``recovery'' with reversible moves increases the number of Nim-heaps in almost-terminal positions. Of course, it also applies to reduced canonical trees more complicated than Nim-heaps: because of reducible moves, the reduced canonical tree of a position can effectively be less complicated that the reduced canonical trees of its children. For example, this is the case for the tree of figure \ref{reversible}, corresponding to the Sprouts position \verb/0.0.A.}1A.}]!/.

\subsection{Factoring by $\nim{1}$}
\label{facto1}

Some reduced canonical trees can be written in the form: $\mathscr{G}+\nim{1}$. It becomes interesting when we consider sums of this kind of tree, because we can use the property $\nim{1}+\nim{1}\sim\nim{0}$ to reduce the size of the trees.

For example, by using that $\nim{3}\sim\nim{2}+\nim{1}$ and that $\{\nim{3};\{\nim{2}\}\}\sim\nim{1}+\{\nim{2}\}$, we obtain: $\nim{3}+\{\nim{3};\{\nim{2}\}\}\sim\nim{2}+\nim{1}+\nim{1}+\{\nim{2}\}\sim\nim{2}+\{\nim{2}\}$. It means that the sum of two trees of height 3 and 4 has been reduced to the sum of two trees of height 2 and 3.

Some other sums enable to reduce the size of the trees. Conway notes for example that  $\{\nim{0};\{\nim{2}\};\{\nim{3};\{\nim{2}\}\}\}+\nim{2}\sim\{\nim{2}\}$ in \cite{jc01} p. 151. But such sums are much too rare to be useful and in our program we used only: $\nim{1}+\nim{1}\sim\nim{0}$.

\section{Computation of reduced canonical trees}

\subsection{Representation and storage}

\subsubsection{String representation}

The intuitive way of representing RCTs\footnote{In the following of this article, \textit{reduced canonical tree} will be shortened \textit{RCT}.} is to recursively define strings, in which each RCT is represented by the set of its children.

For example, the RCT of $S_4$ would be represented by: $\{\nim{3};\{\nim{1};\nim{2};\{\nim{3};\{\nim{2}\}\}\}\}$\footnote{this is a compact form. With a string representation even for Nim-heaps, we should replace $\nim{0}$ with $\{\}$, $\nim{1}$ with $\{\{\}\}$, $\nim{2}$ with $\{\{\};\{\{\}\}\}$ and $\nim{3}$ with $\{\{\};\{\{\}\};\{\{\};\{\{\}\}\}\}$}. But this string representation quickly becomes inefficient when the RCTs grow. If the same RCT occurs in several places of a bigger one, its representation is stored as many times as it occurs, whereas it is clear that it would suffice to store this information only once. This defect is even more important when we store many RCTs in the database.

\subsubsection{Link representation}
\label{linkrep}

To avoid the problem of redundancy in string representations, we implemented a representation by link, where each RCT has an identifier (a number). An RCT is still represented by the set of its children, but this time, we only store the set of their identifiers instead of their complete representation. This allows us to store only once each RCT in the database, while we can still refer to them several times by using their identifiers.

But during computations using the RCTs, we frequently need some information about a given RCT: its height, and its outcome. Of course, it is possible to compute them recursively, but the search of identifiers in the databases would cost much running time.

Consequently, we choose a representation that contains the main information needed about an RCT, and the identifier is made of three parameters:
\begin{itemize}
 \item the height of the RCT.
 \item a number to distinguish RCTs of same height.
 \item the character ``\verb/W/'' or ``\verb/L/'' according to the mis\`ere outcome of the RCT.
\end{itemize}

By convention, the number is \verb/0/ if the RCT is a Nim-heap. Otherwise, the number \verb/1/ is given to the first RCT of a given height that we meet, \verb/2/ to the second one...

\begin{figure}[ht]
\centering
\includegraphics[scale=0.4]{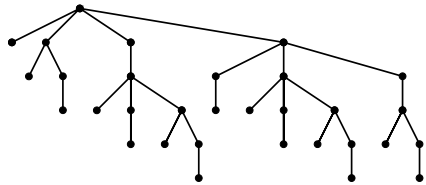}
\caption{Reduced canonical tree of height 5}
\label{exrct}
\end{figure}

The character that describes the outcome of the RCT is useful during the reduction process: a tree is reducible to $\nim{0}$ only if it is a Win in mis\`ere version (cf paragraph \ref{red0}).

let us give an example with the RCT of figure \ref{exrct}, which arises from Sprouts position: \verb/1ABC.}BCDE.}ADE.}]!/. The string representation of this RCT is:\\$\{\nim{0};\nim{2};\{\nim{3}\};\{\nim{1};\nim{3};\{\nim{2}\}\}\}$.

In the following table, for each subtree of this RCT, we give its identifier as well as the set of identifiers of its children.

\medskip

\begin{center}
\begin{tabular}{ |*{3}{c|} }
\hline
RCT & identifier & set of children \tabularnewline
\hline
$\nim{0}$ & \verb/0-0-W/ & - \tabularnewline
\hline
$\nim{1}$ & \verb/1-0-L/ & \verb/0-0-W/ \tabularnewline
\hline
$\nim{2}$ & \verb/2-0-W/ & \verb/0-0-W 1-0-L/ \tabularnewline
\hline
$\nim{3}$ & \verb/3-0-W/ & \verb/0-0-W 1-0-L 2-0-W/ \tabularnewline
\hline
$\{\nim{2}\}$ & \verb/3-1-L/ & \verb/2-0-W/ \tabularnewline
\hline
$\{\nim{3}\}$ & \verb/4-1-L/ & \verb/3-0-W/ \tabularnewline
\hline
$\{\nim{1};\nim{3};\{\nim{2}\}\}$ & \verb/4-2-W/ & \verb/1-0-L 3-0-W 3-1-L/ \tabularnewline
\hline
$\{\nim{0};\nim{2};\{\nim{3}\};\{\nim{1};\nim{3};\{\nim{2}\}\}\}$ & \verb/5-1-W/ & \verb/0-0-W 2-0-W 4-1-L 4-2-W/ \tabularnewline
\hline
\end{tabular}
\end{center}

\medskip

The RCTs met during the computation are stored in a database similar to the last two columns of this table, as ( RCT ; set of children of the RCT ) couples.

\subsubsection{Dependence on the order of computation}

However, this link representation presents a disadvantage: the number in the identifier depends only on the order in which the RCTs have been met, so the same RCT can have different identifiers in different computations. As a consequence, RCT databases produced in different computations are incompatible.

 We explain in paragraph \ref{criterion} that in our computations, we produced an RCT database once and for all to circumvent this problem.

\subsection{Computation of the reduced canonical tree obtained from a given position}

The definition of paragraph \ref{rct} provides immediately a recursive algorithm to compute the RCT of a given position:

\begin{itemize}
\item compute (recursively) the RCT of each child of the position.
\item delete the duplicates amongst these RCTs.
\item reduce the obtained tree, using the smallest possible reducer.
\end{itemize}

There is no difficulty in the first two steps, but the programming of the reduction should be detailed.

\subsubsection{Reduction in the case of Nim-heaps}
\label{rednim}

The first reduction to be considered is the one corresponding to the theorem \ref{theomex}: if all the RCTs obtained from the children of the position are Nim-heaps, and if at least one is $\nim{0}$ or $\nim{1}$, then the RCT of the position is itself a Nim-heap that we can determine with the $Mex$ rule.

\subsubsection{Reduction to $\nim{0}$}
\label{red0}

Next, we test whether the RCT of the position is reducible to $\nim{0}$. We have seen in paragraph \ref{emptytree} that it is possible only if the position is a Win in mis\`ere version.

We must therefore start by checking if a child is a Loss in mis\`ere version, which is immediate, because the results are stored in the identifiers (without this, we should determine the outcome of the position through a computation on the entire tree, which would be much more expensive in running time).

Then, we check if every child is a reversible move, ie if it has $\nim{0}$ as a child.

\subsubsection{Other reducers}

If none of the previous reductions have worked, we have to see if it is possible to find another reducer. So, imagine that we are computing the RCT of a position, and that after having removed the duplicates amongst the RCTs of its children, we have obtained the set: $\{\mathscr{A}_1;\mathscr{A}_2;\mathscr{A}_3...\}$. Let $h_i$ denote the height of $\mathscr{A}_i$. We sort the RCTs $\mathscr{A}_i$ so that $(h_i)$ is an increasing sequence. Then the observation of paragraph \ref{reducers} implies that it suffices to test the reducers of the form $\{\mathscr{A}_1;\mathscr{A}_2;...;\mathscr{A}_i\}$, where $h_{i+1}\geq h_i+2$.

Coming back to the example of paragraph \ref{linkrep}, let us imagine that we want to reduce: \{\verb/0-0-W 2-0-W 4-1-L 4-2-W 5-0-W 7-0-W/\}. We then only need to test the following potential reducers:

\begin{itemize}
\item \{\verb/0-0-W/\}
\item \{\verb/0-0-W 2-0-W/\}
\item \{\verb/0-0-W 2-0-W 4-1-L 4-2-W 5-0-W/\}
\end{itemize}

To be sure of finding the smallest possible reducer, our algorithm examines first the smallest potential reducers.

\{\verb/0-0-W/\}=\verb/1-0-L/ is not a correct reducer, because even if \verb/1-0-L/ is a child of \verb/2-0-W/, it is not a child of \verb/4-1-L/.

Then, our algorithm examines the potential reducer \{\verb/0-0-W 2-0-W/\}. It starts by looking in the database for the identifier of the RCT that has these two children, but it does not find it. For a good reason: \{\verb/0-0-W 2-0-W/\} is reducible to \verb/1-0-L/ (see  paragraph \ref{rednim}). As the smallest possible reducer must be an RCT, and that this RCT must be the son of at least one child of the position, the recursive nature of the algorithm implies that this RCT has already been stored in the database. So, when our algorithm cannot find a potential reducer in the database, it means that this potential reducer is itself reducible, and there is no need to try it.

It remains to test the potential reducer \{\verb/0-0-W 2-0-W 4-1-L 4-2-W 5-0-W/\}, which is not suitable, since the only children of \verb/7-0-W/ are Nim-heaps. So it is impossible to reduce \{\verb/0-0-W 2-0-W 4-1-L 4-2-W 5-0-W 7-0-W/\}. Since it is a new RCT, our program generates a new identifier of the form \verb/8-n-W/ (this new RCT is a Win in mis\`ere version because the child \verb/4-1-L/ is a Loss).

\subsection{Factoring by $\nim{1}$}

We have seen in paragraph \ref{facto1} that it is useful to determine which RCTs could be written as a sum of another RCT and of the Nim-heap $\nim{1}$. We present here the implementation of this factorization. The factorization step must take place just after the various reductions in the recursive algorithm.

\subsubsection{Representation of a position factorizable by $\nim{1}$}

If two RCTs $\mathscr{G}$ and $\mathscr{H}$ conform to $\mathscr{G}\sim\mathscr{H}+\nim{1}$, then the difference between their heights is 1. In order to store this relation in our database, we express the identifier of the biggest RCT as function of the smallest one, with the following notation.

We know that $\nim{3}\sim\nim{2}+\nim{1}$ (or, symetrically, that $\nim{2}\sim\nim{3}+\nim{1}$). The smallest of these two trees is $\nim{2}$. Consequently, $\nim{2}$ keeps the same identifier, \verb/2-0-W/, while the identifier of $\nim{3}$ will be \verb/2-0+1-W/.

We have also seen that $\{\nim{3};\{\nim{2}\}\}\sim\nim{1}+\{\nim{2}\}$. As the identifier of $\{\nim{2}\}$ is \verb/3-1-L/, the identifier of $\{\nim{3};\{\nim{2}\}\}$ will be \verb/3-1+1-W/ (the outcome character ``\verb/W/'' corresponds to \verb/3-1+1/, not to \verb/3-1/).

\subsubsection{Detection of a position factorizable by $\nim{1}$}

We come back to the example of paragraph \ref{linkrep}. Using the new identifier described in the previous paragraph, we get: \verb/4-2-W/$=$\verb/{0-0+1-L 2-0+1-W 3-1-L}/. Now, let us compute the children of \verb/4-2-W/+$\nim{1}$. To get a child, either we play a move in \verb/4-2-W/ or in $\nim{1}$. So we have: \verb/4-2-W/+$\nim{1}$=\verb/{0-0-W 2-0-W 3-1+1-W 4-2-W}/.

In this example, we can observe a more general result: if $\mathscr{G}$ is an RCT that cannot be factorized by $\nim{1}$, then $\mathscr{G}+\nim{1}=\{$(child of $\mathscr{G}$)$+\nim{1}$ ; $\mathscr{G}\}$. So, to determine whether the set of the RCTs of the children of a position corresponds to a position factorizable by $\nim{1}$, it suffices to check whether the previous set is of the correct form. In particular, $\mathscr{G}$ is the only RCT of maximal height that cannot be factorized by $\nim{1}$.

\section{Mis\`ere computation algorithm using RCTs}

When we compute the outcome of a position with our program, we explore the game tree obtained from this position in order to determine the outcome of the root. Details on how we carry out this exploration are available in \cite{lv07}. What changes compared to the normal version is the nature of the nodes of the explored tree.

\subsection{Simplification of the positions with RCTs}
\label{simplif}

We detail the nature of these nodes on an example. Consider the Sprouts position:\\
\verb/0.0.0.0.0.0.0.0.}]22.}]2ab2ba.}]0.0.A.}2A.}]!/

This position has 4 independent positions, more or less complex:

\begin{itemize}
 \item \verb/0.0.0.0.0.0.0.0.}]/ is too big to compute its RCT.
 \item \verb/22.}]/$\sim\nim{1}$.
 \item \verb/2ab2ba.}]/$\sim\nim{3}\sim\nim{1}+\nim{2}\sim\nim{1}+$\verb/2-0-W/.
 \item \verb/0.0.A.}2A.}]/$\sim\nim{1}+$\verb/3-1-L/.
\end{itemize}

Finally, using $\nim{1}+\nim{1}+\nim{1}\sim\nim{1}$, we obtain that this position is indistinguishable of \verb/0.0.0.0.0.0.0.0.}]!/$+\nim{1}+$\{\verb/2-0-W/$+$\verb/3-1-L/\}.

Broadly speaking, a \textit{node} of the game tree is composed of three parts:

\begin{itemize}
 \item the position part, which includes one or more independent positions, too big to compute their RCT.
 \item the $\nim{0}/\nim{1}$ part, which is $\nim{0}$ or $\nim{1}$ according to the parity of the number of $\nim{1}$.
 \item the RCT part, which contains a list of RCTs.
\end{itemize}

\subsection{Computation of the children of a node}

The children of a node are computed by taking into account 3 types of possible moves. A move is done in one of the node part, leaving the other two unaltered :

\begin{itemize}
 \item a normal Sprouts move in the position part
 \item a move in the $\nim{0}/\nim{1}$ part, which consists to replace $\nim{1}$ by $\nim{0}$.
 \item a move in the RCT part, which consists to replace an RCT by one of its children
\end{itemize}

\subsection{Interest of RCTs}

Replacing independent positions by their RCT when they are known has several advantages. Many almost-terminal positions have the same RCT, so replacing them by their RCT can simplify the game tree, by reducing the number of nodes stored and explored (and therefore, it reduces memory consumption and improves the running time).

For bigger positions, we found that RCTs provide relatively few simplifications. First, it is rather rare that complex Sprouts positions have the same RCT, and second, our canonization of Sprouts positions is already quite good: there are very few duplicates amongst the computed children of a given position, so the RCT of this position has more or less the same number of children as this position. But even so, the RCT proves its usefulness, because computing the children of a Sprouts position is complex and expensive in running time, and so the database of ( RCT ; set of children of the RCT ) couples provides a cache that avoids making unnecessarily the same computations of children several times.

\subsection{Criterion to replace a position by its RCT}
\label{criterion}

Unfortunately, computing the RCT of a position requires the exploration of its whole game tree. This limits the size of the positions that can be replaced by their RCT. Indeed, the size of the game tree obtained from a position increases rapidly with the size of the position. It can be seen in the following table where we counted the number of canonical trees obtained from different starting positions of Sprouts.

\begin{figure}[ht]
\centering
\begin{tabular}{ |*{2}{c|} }
\hline
starting spots & number of canonical trees \tabularnewline
\hline
2 & 10 \tabularnewline
\hline
3 & 55 \tabularnewline
\hline
4 & 713 \tabularnewline
\hline
5 & 10461 \tabularnewline
\hline
6 & 150147 \tabularnewline
\hline
\end{tabular}
\caption{Number of different canonical trees in game trees obtained from a starting position}
\label{nbrct}
\end{figure}

We considered two different criteria to decide when to compute the RCT of a position. The first is to compute the RCT of all positions below a certain limit number of lives. This method has the advantage to adapt to the current computation, because we compute only the RCTs of positions actually encountered. But it has the disadvantage of changing the RCT databases during the computation. Above all, this criterion is not well adapted to the structure of Sprouts: two positions with the same number of lives can lead to trees of completely different complexity. For example, there are 55 different canonical trees in the game tree obtained from $S_3$ (ie 9 lives), whereas there are 478 in the tree resulting from the position \verb/1abcde2edcba.2.}]!/ (which also has 9 lives).

We have therefore chosen another criterion. The starting positions with $p$ spots, and their descendants, quickly appear in the computation of positions with a higher number of spots. We therefore computed the RCT of $S_6$. Once this RCT computed, we do not compute any other, and we thus have a fixed RCTs database.

For example, from $S_{12}^-$, we play the move \verb/0.0.0.0.0.0.0.0.AB.}0.0.0.AB.}]!/, and then the move \verb/0.0.0.0.0.0.0.0.}]0.0.A.}0.A.}]!/. At this time, our algorithm modifies the node, because we can read in the database obtained from $S_6$ that \verb/0.0.A.}0.A.}]!/$\sim$\verb/3-1+1-W/. On the other side, the rest of the position does not appear in this database. Thus the new node is \verb/0.0.0.0.0.0.0.0.}]!/$+\nim{1}+$\verb/3-1-L/.

\subsection{Sums of positions}

To simplify the nodes, we could consider merging the $\nim{0}/\nim{1}$ part and the list of RCTs in a single RCT, by computing their sum. In the case of paragraph \ref{simplif}, we would replace $\nim{1}+$\{\verb/2-0-W/$+$\verb/3-1-L/\} by a single RCT of height 6, \verb/5-1+1-W/. The interest of this approach is twofold: in addition to writing simpler nodes, it would detect the simplifications described at the end of paragraph \ref{facto1}.

However, this method is not practical, because the computation of the sum of RCTs requires the storage of too many additional RCTs. For example, after the storage of the RCT of $S_6$, let us suppose that we want to calculate the result of the position \verb/0.0.0.0.0.}]0.0.0.0.}]!/. The RCTs of these two independent positions are known (and their respective heights are 13 and 6).

The algorithm described above would replace each of these two positions by its RCT, and thus start the computation on the node: \verb/!/$+\nim{0}+$\{\verb/13-7-W+6-208-L/\}. Then, the classical Win/Loss algorithm we used allows us to find the result of this position by storing only 10 losing positions, while in contrast, the computation of the RCT of \verb/13-7-W+6-208-L/ needs to store more than 35,000 new RCTs, ie more than for the computation of the RCT of $S_6$.

During a mis\`ere computation, we frequently need to consider such sums (and even more complex ones), so it is not possible to compute the resulting RCTs. Finally, the simplifications expected at the end of paragraph \ref{facto1} do not compensate for this problem: during the tests we conducted, we have observed none.

\section{Results}

\subsection{Phase of preliminary computations}

The following table summarizes the number of RCTs involved in game trees obtained from starting positions, and allows us to imagine the memory needed to take $S_7$ as the basis of our computations instead of $S_6$ (we have not computed it due to memory limitations, but it may be accessible: only 45 MB of RAM are needed to store $S_6$).

\medskip

\begin{center}
\begin{tabular}{ |*{3}{c|} }
\hline
starting spots & number of RCTs & number of positions stored \tabularnewline
\hline
2 & 5 & 18\tabularnewline
\hline
3 & 7 & 157\tabularnewline
\hline
4 & 35 & 1796\tabularnewline
\hline
5 & 1204 & 24784\tabularnewline
\hline
6 & 25459 & 393103\tabularnewline
\hline
\end{tabular}
\end{center}

\medskip

These numbers of RCTs are an objective data of the Sprouts game, which should be checkable by other programmers. In contrast, the number of positions encountered in the game trees obtained from starting positions depends on our Sprouts positions canonization. This number is an indicator of the quality of our canonization. As our canonization only uses simplifications that preserve the canonical trees, it is wise to compare this column with table \ref{nbrct}.

While many almost-terminal positions are indistinguishable from Nim-heaps, some quite complex positions are indistinguishable from Nim-heaps as well. For example, \verb/0.0.0.0.2.}]!/$\sim\nim{0}$. Conversely, the position with a minimum number of lives that is not a Nim-heap is \verb/ABC.}ABD.}CE.}DE.}]!/$\sim\{\nim{2}\}$.

\subsection{Computation of $S_{17}^-$}

The implementation of the techniques described in this article allowed us to compute that $S_{17}^-$ is a Win. The computation required the storage of approximately 170,000 nodes. Of these nodes, about half had an empty position part, which shows that we can explore the game tree so that the position is dismantled rather quickly into several independent positions small enough for their RCT to be in the game tree of $S_6$.

The exploration of the game tree was performed in the same way as what we described in \cite{lv07}: a depth-first algorithm, in which we track the exploration, and manually choose to explore the branches that seem the most hopeful. We also used a check algorithm : when a computation is finished, this algorithm only keeps the positions needed to demonstrate the outcome. It reduces significantly the sizes of the databases: for example, amongst the 170,000 positions computed, less than 18,000 are necessary to demonstrate $S_{17}^-$.

The computation time was only twenty hours on a 1.8 GHz processor, and RAM consumption less than 100 MB, so we can reasonably expect to improve this record. Figure \ref{nbpos} shows the number of positions needed to demonstrate $S_p^-$ (after having used a check algorithm). This number gives a good idea of the complexity of $S_p^-$, so that we could imagine the difficulty of computing higher values of $p$.

\begin{figure}[h]
\centering
\begin{tabular}{ |*{12}{c|} }
\hline
$p$ & 7 & 8 & 9 & 10 & 11 & 12 & 13 & 14 & 15 & 16 & 17\tabularnewline
\hline
Outcome of $S_p^-$ & L & L & L & W & W & W & L & L & L & W & W\tabularnewline
\hline
Positions & 7 & 24 & 44 & 114 & 79 & 983 & 1082 & 3517 & 6906 & 8981 & 17583\tabularnewline
\hline
\end{tabular}
\caption{Number of positions needed to demonstrate $S_p^-$.}
\label{nbpos}
\end{figure}

It is worthwhile to be aware that this number is the number of positions we need \textit{in addition} to the game tree of $S_6$ to compute the outcome of $S_p^-$. Of course, we can't find the outcome of $S_7^-$ with only 7 positions...

The programming methods described in this article will probably help to determine the outcome of $S_{18}^-$ or $S_{19}^-$ fast enough, but it is unlikely that the study of the mis\`ere game will give results equivalent to the normal version. Indeed, because of theoretical difficulties related to the mis\`ere version, many more nodes must be explored to demonstrate a result. Also, unlike in the case of normal Sprouts game, in which the complexity of $S_p$ does not increase strictly with $p$ (e.g. it is faster to compute $S_{17}^+$ than $S_{15}^+$, and we were able to compute $S_{47}^+$ but not $S_{33}^+$), the difficulty seems to increase more regularly for the mis\`ere version.

\section{Conclusion}

The theory of reduced canonical trees and the algorithm that results for mis\`ere impartial games were particularly effective in the case of Sprouts, by using efficiently the divisions in sums of independent positions. The same algorithm could be used for computing other impartial games in mis\`ere version, provided that such divisions in sums of independent position occur in these games.

The potential improvements may relate to various areas: in addition to purely theoretical improvements, and of course to improvements due to the increasing performance of computers, we may hope for improvements in the game trees exploration.

The program we used for computation is available with its source code on our web site \verb|http://sprouts.tuxfamily.org/| under a GNU licence, together with several databases used during our computations.

\bibliography{misere_en}

\providecommand{\bysame}{\leavevmode\hbox to3em{\hrulefill}\thinspace}
\providecommand{\MR}{\relax\ifhmode\unskip\space\fi MR }
% \MRhref is called by the amsart/book/proc definition of \MR.
\providecommand{\MRhref}[2]{%
  \href{http://www.ams.org/mathscinet-getitem?mr=#1}{#2}
}
\providecommand{\href}[2]{#2}
\begin{thebibliography}{1}

\bibitem{ajs91}
D.~Applegate, G.~Jacobson, and D.~Sleator, \emph{Computer {A}nalysis of
  {S}prouts}, Tech. Report CMU-CS-91-144, Carnegie Mellon University Computer
  Science Technical Report, 1991.

\bibitem{ww01}
Edwyn Berkelamp, John Conway, and Richard Guy, \emph{Winning ways for your
  mathematical plays}, A K Peters, 2001.

\bibitem{jc01}
John~H. Conway, \emph{On numbers and games (second edition)}, A K Peters, 2001.

\bibitem{mg67}
Martin Gardner, \emph{Mathematical games : of sprouts and brussels sprouts,
  games with a topological flavor}, Scientific American \textbf{217} (July
  1967), 112--115.

\bibitem{gs56}
P.M. Grundy and C.A.B. Smith, \emph{Disjunctive games with the last player
  losing}, Mathematical Proceedings of the Cambridge Philosophical Society
  \textbf{52} (1956), no.~3, 527--533.

\bibitem{lv07}
Julien Lemoine and Simon Viennot, \emph{A further computer analysis of
  sprouts},
  \url{http://download.tuxfamily.org/sprouts/sprouts-lemoine-viennot-070407.pd%
f}, 2007.

\bibitem{tp05}
Thane~E. Plambeck, \emph{Taming the wild in impartial combinatorial games},
  INTEGERS: Electronic Journal of Combinatorial Number Theory \textbf{5}
  (2005), G5.

\bibitem{si08}
Aaron~N. Siegel, \emph{Mis\`ere games and mis\`ere quotients},  (2006),
  \url{http://arxiv.org/abs/math.CO/0612616}.

\end{thebibliography}
\bibliographystyle{amsplain}

\end{document}